\documentclass[11pt]{amsart}

\usepackage{amsmath,amssymb,url}
\usepackage[mathscr]{eucal}
\usepackage{pstricks,pst-plot}
\usepackage{fullpage}
\newtheorem{thm}{Theorem}

\newrgbcolor{darkgreen}{0.15 0.5 0}
\newrgbcolor{cyan}{0.15 0.2 0.8}
\newcommand{\sout}[1]{}

\newcommand{\Gammaf}{\mathrm{\Gamma}}

\newcommand{\dstyle}{\displaystyle}

\renewcommand{\appendix}{%
\renewcommand{\thesection}{\Alph{section}}%
\section}

\title{Factoring with Hints}

\author{Francesco Sica}
\address{School of Science and Technology\\
Nazarbayev University\\
53 Kabanbay Batyr\\
Astana, 010000 Kazakhstan
}
\email{francesco.sica@nu.edu.kz}
\thanks{Research supported in part by a grant from the Social Development Fund.}

\date{\today}


\begin{document}

\begin{abstract}

We introduce a new deterministic factoring algorithm, which could be described in the cryptographically fashionable term of ``factoring with hints": we show that, given the knowledge of the factorisations of $O(N^{1/3+\epsilon})$ terms surrounding $N=pq$ product of two large primes, we can recover deterministically $p$ and $q$ in $O(N^{1/3+\epsilon})$ bit operations. Although this is slower than the current best factoring algorithms, this method shows that the factorisations of close integers are related and that consequently one can expect more results along this line of thought.

\end{abstract}

\maketitle

\section{Introduction}

The problem of quickly factoring large integers is central in cryptography and computational number theory.
%
The current state of the art in factoring large integers $N$, the Number Field Sieve algorithm~\cite{NFS, NFSdev}, stems from the earlier Quadratic Sieve~\cite{QS} and Continued Fraction~\cite{MorBri75}. We should also mention the Elliptic Curve Method by H.\@ Lenstra~\cite{ECM}, which is particularly useful when $N$ has a small prime factor $p$. They are all probabilistic factoring algorithms.

These algorithms have \emph{heuristic} running times respectively
$O\bigl(\exp(c(\log N)^{1/3} (\log\log N)^{2/3})\bigr)$, $O\bigl(\exp(c(\log N)^{1/2} (\log\log N)^{1/2})\bigr)$ and
$O\bigl(\exp(c(\log p)^{1/2} (\log\log p)^{1/2})\bigr)$, for some constant $c$ (not always the same).
The first two strive to find nontrivial arithmetical relations of the form $x^2\equiv y^2 \pmod N$ (which lead to a nontrivial factor by computing $\gcd(N,x+y)$), whereas the third is a generalisation of Pollard's $p-1$ method~\cite{Polp-1}, involving computations in some elliptic curve group instead of $\mathbb{Z}/N$. We should note, however, that there exist probabilistic algorithms with proved running time $O\bigl(\exp((1+o(1))(\log N)^{1/2} (\log\log N)^{1/2})\bigr)$~\cite{LePo92}. As far as the author is aware, no such rigorous bound exists in the form $O\bigl(\exp\left((\log N)^c\right)\bigr)$ for $c<1/2$. Similarly, no deterministic subexponential algorithm is currently known, the best one being Shank's square form factorization SQUFOF which runs in $O(N^{1/4+\epsilon})$, or in $O(N^{1/5+\epsilon})$ on the Extended Riemann Hypothesis.

In this work, we want to introduce a new paradigm in integer factorisation, one that doesn't supersede previous efforts, but rather complements it by showing that the factorisation of a small number of consecutive integers in related in a nontrivial way. Therefore, if numbers close to a product $N=pq$ of two primes are easier to factor than $N$ itself, we can expect a reduction in the time to factor $N$. Quantifying the number of consecutive integers versus the additional computational effort to find the needed relation is a matter of further investigation, some of which will come out in another work by ourselves. Here we content ourselves with a first nontrivial result.

\begin{thm}
Let $N=pq$ a product of two primes. Then, given an arbitrary $\epsilon>0$, the factors $p$ and $q$ can be recovered in $O(N^{1/3+\epsilon})$ bit operations from the knowledge of the factorisations of $O(N^{1/3+\epsilon})$ integers closest to $N$.
\end{thm}

\section{Notations}
This work borrows heavily from standard notations in analytic number theory and indeed a classical reference on the subject is the treatise of Davenport~\cite{Da}. In particular, we will make liberal use of the $O$ notation in Landau's as well as Vinogradov's form ($\ll$). Hence, for instance
$$
f(u) = O\bigl(g(u)\bigr) \iff f(u)\ll g(u)
$$
means that $g(u)>0$ and $|f(u)|/g(u)$ is bounded above (usually as $u\to\infty$ or $u\to0^+$, depending clearly on the context). Similarly, $f(u)= o(g(u))$ (resp. $f(u)= \Omega(g(u))$) means $g(u)>0$ and $|f(u)|/g(u)$ goes to zero (resp. $|f(u)|/g(u)$ is bounded below).
Unless specified, the implied constants are absolute.

Any sum such as
$$
\sum_{abc=n} a^2bc
$$
is to be understood as taken over all positive integers $a,b,c$ such that $abc=n$. We also define
$$
\sum_{a\mid n} f(a) = \sum_{ab=n} f(a)
$$
so that for instance the number of divisors of $n$ is $\sum_{d\mid n} 1$ and its sum of divisors $\sum_{d\mid n} d$. We also write $s=\sigma + it$, with $\sigma,t\in\mathbb{R}$, according to the established convention in analytic number theory.

We will also write $f^{(n)}$ for the $n$-th derivative of the function $f$.

Finally, we write $a\doteq b$ to signify that $a=b+$terms that are not necessarily negligible in size but can be computed in polynomial time (in the bit size of the challenge to be factored), so that they are negligible in time.

\section{Choice of a Multiplicative Function}
\label{S:multfct}
For $\lambda\in\mathbb R$ define
$$
\sigma_{\lambda}(n)  =
\sum_{d\mid n} d^{\lambda} \enspace.
$$
Our goal will be to compute $\sigma_{1/2}(N) = 1 + \sqrt N + \sqrt p + \frac{\sqrt N}{\sqrt p}$ within $O(1/N)$. If so, then one gets an approximation $\mathcal A$ to
\begin{equation}
\label{E:approx}
f(p) = \sqrt p + \frac{\sqrt N}{\sqrt p} = \mathcal A + O\left(\frac1N\right) \enspace.
\end{equation}
Let us study the function in $(0, \infty)$
$$
f(z) = \sqrt z + \frac{\sqrt N}{\sqrt z} \Rightarrow f'(z) = \frac{1}{2\sqrt z} \left( 1 - \frac{\sqrt N}{z}\right)
\Rightarrow f''(z) = \frac{1}{4z^{3/2}}\left( \frac{3\sqrt N}{z} - 1\right)
\enspace.
$$
The function $f$ is convex in $(0, 3\sqrt N)$ with a unique critical point (and therefore absolute minimum) at $z= \sqrt N$. We will suppose that $N=pq$ with $p<\sqrt N < q$. In fact, we may as well suppose that $p\leq \sqrt N-2$ by inspection.  Note that $f''(z) \geq 2 N^{-3/4}$ for $z\leq \sqrt N$ and therefore $|f'(z)| \geq 2 N^{-3/4}$ for $z\leq \sqrt N-1$. Define $a\in(0,\sqrt N-1]$ by $f(a)=\mathcal A$.
To see that such $a$ exists, notice that $f$ is decreasing in $(0,\sqrt N]$. If $\mathcal A < f(\sqrt N-1)$, then, for some $\theta\in(\sqrt N-2, \sqrt N-1)$, we can write
$$
\left| \mathcal A - f(p) \right| < | f(\sqrt N-1) - f(\sqrt N-2)| = |f'(\theta)| \geq \frac2{N^{3/4}} \enspace,
$$
contradicting~\eqref{E:approx}. Given then $a\in(0,\sqrt N-1]$ with $f(a)=\mathcal A$, we obtain, for some $\xi\leq \sqrt N-1$,
$$
\left|a-p\right|\, \left|f'(\xi)\right|=\left| f(a) - f(p) \right| = \left| \mathcal A - f(p) \right| \ll \frac1N
\Rightarrow
p=a + O\left(\frac1{N^{1/4}}\right)
$$
and therefore $p=\lfloor a \rceil$, the integer nearest to $a$.

\section{Choice of a Test Function}
\label{S:3}

Consider the Riemann zeta function
$$
\zeta(s)= \sum_{n\geq 1} \frac{1}{n^s} \enspace,
$$ convergent for $\Re s>1$. Then

$$
\zeta(s)\zeta(s-1/2) = \sum_{n=1}^\infty \frac{\sigma_{1/2}(n)}{n^s} \enspace,
$$
absolutely convergent whenever $\Re s > 3/2$.
Now let for $\nu\in\mathbb N$\footnote{In fact, $\nu$ doesn't need to be an integer, but it simplifies calculations to assume so.} with $\nu\geq2$,
$$
f(t)= \begin{cases}
(1-t)^{\nu-1} & 0\leq t\leq 1,\\
0 & t\geq 1.
\end{cases}
$$

The Mellin transform of $f$ is by definition the beta function

$$
B(\nu,s)=\frac{\Gammaf(\nu)\Gammaf(s)}{\Gammaf(s+\nu)} =
\int_0^\infty f(t) t^{s-1}\,dt
$$
hence by the inverse Mellin transform\footnote{We will also use the notation $\dstyle \int\limits_{(c)}$ instead of $\dstyle\int_{c-i\infty}^{c+i\infty}$.},

\begin{equation}
\label{E:3} \frac1{2\pi i}\int_{5/2-i\infty}^{5/2+i\infty}\zeta(s)\zeta(s-1/2) \,
\frac{\Gammaf(\nu)\Gammaf(s)}{\Gammaf(s+\nu)}x^s\,ds
= \sum_{n\leq x} \sigma_{1/2}(n) f\left(\frac
nx\right).
\end{equation}

Call the right-hand side
$$
F_\nu(x) = \sum_{n\leq x}\sigma_{1/2}(n) f\left(\frac
nx\right)=\sum_{n\leq x}\sigma_{1/2}(n) \left(1-\frac nx\right)^{\nu-1}
$$
and note that
$$
P_\nu(x)= x^{\nu-1} F_\nu(x)= \sum_{n\leq x} \sigma_{1/2}(n)(x-n)^{\nu-1}
$$
is a piecewise polynomial (given by a different expression between consecutive integers).

\section{Functional Equation of the Riemann Zeta Function}

The Riemann zeta function is a meromorphic function having a simple pole with residue $1$ at $s=1$ and satisfying the functional equation (given here in asymmetric form)

$$
\zeta(s)= \frac1{2\pi i} \left({2\pi}\right)^s \, \Gammaf(1-s) \zeta(1-s)
\left(e^{i\pi s/2} - e^{-i\pi s/2}\right)
= \frac1{\pi} \left({2\pi}\right)^s \, \Gammaf(1-s) \zeta(1-s) \sin \frac{\pi s}2
\enspace.
$$

\section{Another Expression for $P_\nu(x)$}

A standard ``integration line moving"  (to $\Re s=-1/4$) argument in the integral of~\eqref{E:3} will get us to the following, after picking up the residues of the integrand at $s=3/2$, $s=1$ and $s=0$,

$$
F_\nu(x) = \zeta(3/2)\frac{\Gammaf(\nu)\Gammaf(3/2)}{\Gammaf(\nu+3/2)}x^{3/2}
+ \frac{\zeta(1/2)}\nu x
+ \zeta(0)\zeta(-1/2) +
\frac1{2\pi i} \int\limits_{(-1/4)} \zeta(s)\zeta(s-1/2) \,
\frac{\Gammaf(\nu)\Gammaf(s)}{\Gammaf(s+\nu)}x^s\,ds \enspace.
$$

In fact, we can move the line of integration to $\Re s=-1/4$, since to the right of that line, for any given $\epsilon>0$,
$$
|\zeta(s)\zeta(s-1/2)| \ll |t|^{2+\epsilon} \enspace,
$$
while
$$
\left|\frac{\Gammaf(s)}{\Gammaf(s+\nu)}\right| \ll |t|^{-\nu} \enspace.
$$
In particular the integral on the right-hand side is absolutely convergent when $\nu\geq 4$.
It is this integral is the next focus of our investigation. It is very natural at this point to use the functional equation. We get quite straightforwardly

\begin{align*}
&\frac1{2\pi i} \int\limits_{(-1/4)} \zeta(s)\zeta(s-1/2) \,
\frac{\Gammaf(\nu)\Gammaf(s)}{\Gammaf(s+\nu)}x^s\,ds \\
&=
\frac{1}{(2\pi i)^3} \int\limits_{(-1/4)} \left({2\pi}\right)^s \Gammaf(1-s) \zeta(1-s) \sin\frac{\pi s}2
\\
&\times
\left({2\pi}\right)^{s-1/2} \Gammaf(3/2-s) \zeta(3/2-s)  \sin\left(\frac{\pi s}2 -\frac\pi4\right)
\frac{\Gammaf(\nu)\Gammaf(s)}{\Gammaf(s+\nu)}x^s\,ds \\
&=
\frac{e^{-i\pi /4}}{(2\pi i)^3\sqrt{2\pi}} \int\limits_{(-1/4)} ({4\pi^2 x})^s e^{i\pi s} \Gammaf(1-s) \zeta(1-s)
 \Gammaf(3/2-s) \zeta(3/2-s) \frac{\Gammaf(\nu)\Gammaf(s)}{\Gammaf(s+\nu)}\,ds \\
&+
\frac{e^{i\pi /4}}{(2\pi i)^3\sqrt{2\pi}} \int\limits_{(-1/4)} ({4\pi^2 x})^s e^{-i\pi s} \Gammaf(1-s) \zeta(1-s)
 \Gammaf(3/2-s) \zeta(3/2-s) \frac{\Gammaf(\nu)\Gammaf(s)}{\Gammaf(s+\nu)}\,ds \\
&-
\frac{1}{(2\pi i)^3\sqrt{\pi}}
\int\limits_{(-1/4)} ({4\pi^2 x})^s  \Gammaf(1-s) \zeta(1-s)
 \Gammaf(3/2-s) \zeta(3/2-s) \frac{\Gammaf(\nu)\Gammaf(s)}{\Gammaf(s+\nu)}\,ds \\
&=
\frac{ e^{-i\pi /4}x}{2\pi i \sqrt{2\pi}}  \int\limits_{(5/4)} ({4\pi^2 x})^{-s} e^{-i\pi s} \Gammaf(s) \zeta(s)
 \Gammaf(s+1/2) \zeta(s+1/2) \frac{\Gammaf(\nu)\Gammaf(1-s)}{\Gammaf(1-s+\nu)}\,ds \\
&+
\frac{ e^{i\pi /4}x}{2\pi i \sqrt{2\pi}}  \int\limits_{(5/4)} ({4\pi^2 x})^{-s} e^{i\pi s} \Gammaf(s) \zeta(s)
 \Gammaf(s+1/2) \zeta(s+1/2) \frac{\Gammaf(\nu)\Gammaf(1-s)}{\Gammaf(1-s+\nu)}\,ds \\
&+
\frac{x}{2\pi i\sqrt\pi} \int\limits_{(5/4)} ({4\pi^2 x})^{-s} \Gammaf(s) \zeta(s)
 \Gammaf(s+1/2) \zeta(s+1/2) \frac{\Gammaf(\nu)\Gammaf(1-s)}{\Gammaf(1-s+\nu)}\,ds \enspace.
\end{align*}
Using the Legendre duplication formula
$$
\Gammaf(s)\Gammaf\left(s+\frac12\right) = \sqrt\pi\, 2^{1-2s} \Gamma(2s) \enspace,
$$
together with the functional equation $s\Gammaf(s)=\Gammaf(s+1)$, we obtain
\begin{equation}
\label{E:Gammastoone}
\Gammaf\left(s+\frac12\right) \frac{\Gammaf(s)\Gammaf(1-s)}{\Gammaf(1-s+\nu)} = \Gammaf\left(s+\frac12\right)\Gammaf(s-\nu)
(\cos \pi\nu - \sin\pi\nu \cot\pi s) =  (-1)^\nu \frac{\sqrt\pi\, 2^{1-2s}\Gammaf(2s)}{(s-1)(s-2)\cdots (s-\nu)}
\enspace.
\end{equation}
We can further transform~\eqref{E:Gammastoone} by noting that there exist unique constants $c_{0,\nu}=1, \dots, c_{\nu,\nu}$ such that
\begin{multline*}
\frac{1}{(s-1)(s-2)\cdots (s-\nu)} = \frac{2^\nu}{(2s-2)(2s-4)\cdots (2s-2\nu)} \\
=
\frac{2^\nu c_{0,\nu}}{(2s-1)(2s-2)\cdots(2s-\nu)} + \frac{2^\nu c_{1,\nu}}{(2s-1)(2s-2)\cdots(2s-(\nu+1))}
+\cdots
+ \frac{2^\nu c_{\nu,\nu}}{(2s-1)(2s-2)\cdots(2s-2\nu)}
\end{multline*}
whenever this expression makes sense. To see this, multiply both sides by $(2s-1)(2s-2)\cdots(2s-2\nu)$. The resulting left-hand side is a polynomial of degree $\nu$, expressed as a linear combination of the polynomials resulting from the right-hand side, which form a basis of the vector space of polynomials of degree $\leq \nu$. For instance, $c_{0,1}= c_{1,1}=1$ and $c_{0,4}=1,c_{1,4}=10, c_{2,4}=45, c_{3,4}=c_{4,4}=105$. From~\eqref{E:Gammastoone} we get
$$
(-1)^\nu \Gammaf\left(s+\frac12\right)\Gammaf(s-\nu) = (-1)^\nu \sqrt\pi\, 2^{\nu+1-2s}
\sum_{m=0}^\nu c_{m,\nu} \Gammaf(2s-\nu-m) \enspace.
$$

Putting it together we obtain
\begin{align*}
&\frac1{2\pi i} \int\limits_{(-1/4)} \zeta(s)\zeta(s-1/2) \,
\frac{\Gammaf(\nu)\Gammaf(s)}{\Gammaf(s+\nu)}x^s\,ds \\
&=(-1)^\nu 2^{\nu+1/2}  e^{-i\pi/4} x \Gammaf(\nu)
\sum_{m=0}^\nu c_{m,\nu} \sum_{n\geq 1} \sigma_{-1/2}(n) \frac1{2\pi i}
\int\limits_{(5/4)} ({16\pi^2 xn})^{-s} e^{-i\pi s} \Gammaf(2s-\nu-m)\,ds
\\
&+(-1)^\nu 2^{\nu+1/2}  e^{i\pi/4} x \Gammaf(\nu)
\sum_{m=0}^\nu c_{m,\nu} \sum_{n\geq 1} \sigma_{-1/2}(n) \frac1{2\pi i}
\int\limits_{(5/4)} ({16\pi^2 xn})^{-s} e^{i\pi s} \Gammaf(2s-\nu-m)\,ds
\\
&+(-1)^\nu 2^{\nu+1}  x \Gammaf(\nu)
\sum_{m=0}^\nu c_{m,\nu} \sum_{n\geq 1} \sigma_{-1/2}(n) \frac1{2\pi i}
\int\limits_{(5/4)} ({16\pi^2 xn})^{-s}  \Gammaf(2s-\nu-m)\,ds \enspace.
\end{align*}

We have, for $y>0$,
$$
\frac1{2\pi i} \int\limits_{(1/2)} y^{-s} \Gammaf(s)\,ds= e^{-y}
$$
and after collecting the residues of the gamma function at the negative integers,
$$
\frac1{2\pi i} \int\limits_{(-(2k+1)/2)} y^{-s} \Gammaf(s)\,ds= e^{-y} - \left(
1-y+\frac{y^2}{2!} + \cdots + (-1)^k \frac{y^k}{k!}
\right) \qquad k\geq 0 \enspace.
$$
Remark that if $k\geq 1$, since $|\Gammaf(s)| < e^{-\pi |t|/2}|t|^{-k-1}$ on $\Re s=-k-1/2$, the left-hand side of the previous expression is analytic for $\Re y>0$ and continuous up to $\Re y=0$. Therefore the previous formula for $k\geq 1$ holds in the closed half-plane $\Re y\geq 0$. With this explicit expression we find that 
\begin{align}
&\frac1{2\pi i} \int\limits_{(-1/4)} \zeta(s)\zeta(s-1/2) \,
\frac{\Gammaf(\nu)\Gammaf(s)}{\Gammaf(s+\nu)}x^s\,ds  \notag
\\
&=(-1)^\nu 2^{\nu+1/2}  e^{-i\pi/4} x \Gammaf(\nu)
\sum_{m=0}^\nu c_{m,\nu} \sum_{n\geq 1} \sigma_{-1/2}(n) \frac1{2\pi i}
\int\limits_{(5/4)} ({16\pi^2 xn})^{-s} e^{-i\pi s} \Gammaf(2s-\nu-m)\,ds \notag
\\
&+(-1)^\nu 2^{\nu+1/2}  e^{i\pi/4} x \Gammaf(\nu)
\sum_{m=0}^\nu c_{m,\nu} \sum_{n\geq 1} \sigma_{-1/2}(n) \frac1{2\pi i}
\int\limits_{(5/4)} ({16\pi^2 xn})^{-s} e^{i\pi s} \Gammaf(2s-\nu-m)\,ds \notag
\\
&+(-1)^\nu 2^{\nu+1}  x \Gammaf(\nu)
\sum_{m=0}^\nu c_{m,\nu} \sum_{n\geq 1} \sigma_{-1/2}(n) \frac1{2\pi i}
\int\limits_{(5/4)} ({16\pi^2 xn})^{-s}  \Gammaf(2s-\nu-m)\,ds \notag
\\
&=(-1)^\nu 2^{\nu+1/2}  e^{-i\pi/4} x \Gammaf(\nu)
\sum_{m=0}^\nu \frac{c_{m,\nu}}{2(4\pi i)^{\nu+m} x^{\nu/2+m/2}} \sum_{n\geq 1} \sigma_{-1/2}(n)
\frac{e^{-4\pi i \sqrt{xn}}}{n^{\nu/2+m/2}} \label{E:SSa}
\\
&- (-1)^\nu 2^{\nu-1/2}  e^{-i\pi/4} x \Gammaf(\nu)
\sum_{m=0}^\nu {c_{m,\nu}} \sum_{k=0}^{\nu-1+m} \frac{(-1)^{\nu-3+m-k}}{(\nu-3+m-k)! (4\pi i)^{3+k} x^{3/2+k/2}}\sum_{n\geq 1}
\frac{\sigma_{-1/2}(n)}{n^{3/2+k/2}}   \notag
\\
&+(-1)^\nu 2^{\nu+1/2}  e^{i\pi/4} x \Gammaf(\nu)
\sum_{m=0}^\nu \frac{(-1)^{\nu+m}c_{m,\nu}}{2(4\pi i)^{\nu+m} x^{\nu/2+m/2}} \sum_{n\geq 1} \sigma_{-1/2}(n)
\frac{e^{4\pi i \sqrt{xn}}}{n^{\nu/2+m/2}} \label{E:SSb}
\\
&- (-1)^\nu 2^{\nu-1/2}  e^{i\pi/4} x \Gammaf(\nu)
\sum_{m=0}^\nu {c_{m,\nu}} \sum_{k=0}^{\nu-1+m} \frac{(-1)^{\nu+m}}{(\nu-3+m-k)! (4\pi i)^{3+k} x^{3/2+k/2}}\sum_{n\geq 1}
\frac{\sigma_{-1/2}(n)}{n^{3/2+k/2}}   \notag
\\
&+(-1)^\nu 2^{\nu+1}   x \Gammaf(\nu)
\sum_{m=0}^\nu \frac{c_{m,\nu}}{2(4\pi )^{\nu+m} x^{\nu/2+m/2}} \sum_{n\geq 1} \sigma_{-1/2}(n)
\frac{e^{-4\pi  \sqrt{xn}}}{n^{\nu/2+m/2}}  \notag
\\
&- (-1)^\nu 2^{\nu}   x \Gammaf(\nu)
\sum_{m=0}^\nu {c_{m,\nu}} \sum_{k=0}^{\nu-1+m} \frac{(-1)^{\nu-3+m-k}}{(\nu-3+m-k)! (4\pi )^{3+k} x^{3/2+k/2}}\sum_{n\geq 1}
\frac{\sigma_{-1/2}(n)}{n^{3/2+k/2}} \notag
 \enspace.
\end{align}
In this last expression, the only terms that we cannot calculate explicitly are the two inner series in~\eqref{E:SSa} and~\eqref{E:SSb}.

\section{Factoring with Hints}

We show here, given $\epsilon>0$, how to calculate in $O(N^{1/3+\epsilon})$ bit operations, assuming the factorisation knowledge of $O(N^{1/3+\epsilon})$ integers immediately around $N=pq$, the quantity $\sigma_{1/2}(N)=\sqrt N + 1 + \sqrt p + \sqrt q$ within $O(N^{-1})$, which is sufficient to derive $p$ and $q$. In the following, we suppose that $\nu$ is a fixed (in terms of $N$) integer with $\nu\geq 2+4/3\epsilon$. The work done in the previous section allows us to write

\begin{align}
P_\nu(x) &\doteq (-1)^\nu 2^{\nu-1/2}  e^{-i\pi/4} x^{\nu/2} \Gammaf(\nu)
\sum_{m=0}^\nu \frac{c_{m,\nu}}{(4\pi i)^{\nu+m} x^{m/2}} \sum_{n\geq 1} \sigma_{-1/2}(n)
\frac{e^{-4\pi i \sqrt{xn}}}{n^{\nu/2+m/2}}  \label{E:ST1}
\\
&+
(-1)^\nu 2^{\nu-1/2}  e^{i\pi/4} x^{\nu/2} \Gammaf(\nu)
\sum_{m=0}^\nu \frac{(-1)^{\nu+m}c_{m,\nu}}{(4\pi i)^{\nu+m} x^{m/2}} \sum_{n\geq 1} \sigma_{-1/2}(n)
\frac{e^{4\pi i \sqrt{xn}}}{n^{\nu/2+m/2}} \enspace. \label{E:ST2}
\end{align}

Having fixed $\epsilon>0$, let us approximate the series
$$
\sum_{n\geq 1} \sigma_{-1/2}(n)
\frac{e^{\pm 4\pi i \sqrt{xn}}}{n^{\nu/2+m/2}}
$$
by its $[x^{1/3+2\epsilon}]$-th partial sum, with corresponding error
$$
\ll x^{-(1/3+2\epsilon)(\nu/2-1)}
$$
and therefore~\eqref{E:ST1} and~\eqref{E:ST2} can be replaced by the corresponding expressions where the inner sums in $n$ are truncated at $n\leq [x^{1/3+2\epsilon}]$ with a total error
$$
\ll_\nu x^{\nu/2} x^{-(1/3+2\epsilon)(\nu/2-1)} \leq x^{\nu/3} \enspace.
$$
Remark that the truncated series
$$
\sum_{n\leq y} \sigma_{-1/2}(n)
\frac{e^{\pm 4\pi i \sqrt{xn}}}{n^{\nu/2+m/2}}
=
\sum_{n_1 n_2 \leq y} \frac{e^{\pm 4\pi i \sqrt{xn_1n_2}}}{n_1^{\nu/2+m/2} n_1^{\nu/2+(m+1)/2} }
$$
can be computed trivially in $O(y^{1+\epsilon})$ bit operations since there are $O(y\log y)$ pairs $n_1,n_2$ with $n_1n_2\leq y$.

{ Let $h>0$ and define $\nabla_h P_\nu(x)(= \nabla_h^1 P_\nu(x)) = P_\nu(x)- P_\nu(x-h)$ and $\nabla_h^{k+1} P_\nu(x) = \nabla_h\nabla_h^k P_\nu(x)$ for $k\geq 1$. The following statements can easily be shown by induction.

\begin{enumerate}
\item If $P$ is a polynomial of degree $d$ then $\nabla_h^{d+1} P =0$,

\item $
\nabla_h^k P_\nu(x) = \sum_{i=0}^k \binom ki P_\nu(x-ih)
$.
\end{enumerate}

Letting $x=N+N^{1/3}$ and $h=N^{1/3}$ we see that
$\nabla_h^\nu P_\nu(x)$ can be expressed as $\sigma_{1/2}(N) N^{(\nu-1)/3} +$ terms involving only $\sigma_{1/2}(n)$ for $x-\nu N^{1/3} \leq n\leq x$. On the other hand from~\eqref{E:ST1} and~\eqref{E:ST2} it can also be computed within $O(x^{\nu/3})$ with $O(x^{1/3+3\epsilon})$ bit operations.
}

%


Therefore we can compute in $O(x^{1/3+3\epsilon})$ bit operations an approximation of $x^{(\nu-1)(1/3+\epsilon)}\sigma_{1/2}(N)$  within $O(x^{\nu/3})$. This leads to an approximation of $\sigma_{1/2}(N)$  within $O(1/N)$ as required. Recovering $p\mid N$ is explained in Section~\ref{S:multfct}.

\section{Final Considerations}
Our method relates the factorisations of $O(N^{\theta+\epsilon})$ numbers close to $N$ with the factorisations of the first $O(N^{\theta+\epsilon})$ integers.
The result given here (with $\theta=1/3$) is rather crude, because the series~\eqref{E:ST1} and~\eqref{E:ST2} were approximated by computing a partial sum trivially.
There are two ways in which one can hope to improve this approach: one would be to find a computationally less expensive summation method for these series, while another would be to group together in the partial sum the summation of nearby terms into an expression easy to evaluate.


\bibliography{biblio}
\bibliographystyle{plain}

\end{document}